\documentclass[a4paper,11pt,twoside]{amsart}

\usepackage[T1]{fontenc}
\usepackage[latin1]{inputenc}
\usepackage{xspace} 
\usepackage{amsfonts}
\usepackage{amsmath,amssymb,amsthm,amscd}
\usepackage[all]{xy}

\title{Proper actions on corank-one reductive homogeneous spaces}
\author{Fanny Kassel}
\address{D\'epartement de Math\'ematiques,
B\^atiment~425,
Facult\'e des Sciences d'Orsay,
Universit\'e Paris-Sud 11,
91405 Orsay Cedex,
France}
\email{fanny.kassel@math.u-psud.fr}

\theoremstyle{plain}
\newtheorem{prop}{Proposition}[section]
\newtheorem{theo}[prop]{Theorem}
\newtheorem{coro}[prop]{Corollary}
\newtheorem{lem}[prop]{Lemma}
\theoremstyle{definition}

\newcommand{\C}{\mathbb{C}}
\newcommand{\R}{\mathbb{R}}
\newcommand{\Q}{\mathbb{Q}}
\newcommand{\Z}{\mathbb{Z}}
\newcommand{\N}{\mathbb{N}}
\newcommand{\F}{\mathbb{F}}
\newcommand{\kkk}{\mathbf{k}}

\newcommand{\SLn}{\mathrm{SL}_n}
\newcommand{\SL}{\mathrm{SL}}
\newcommand{\SO}{\mathrm{SO}}

\newcommand{\U}{\mathrm{U}}
\newcommand{\SU}{\mathrm{SU}}

\newcommand{\Sp}{\mathrm{Sp}}
\newcommand{\GLnC}{\mathrm{GL}_n(\mathbb{C})}
\newcommand{\GLnR}{\mathrm{GL}_n(\mathbb{R})}

\newcommand{\PSL}{\mathrm{PSL}}
\newcommand{\M}{\mathrm{M}}
\newcommand{\g}{\mathfrak{g}}
\newcommand{\kk}{\mathfrak{k}}
\newcommand{\p}{\mathfrak{p}}
\newcommand{\aaa}{\mathfrak{a}}
\DeclareMathOperator{\Res}{Res}
\DeclareMathOperator{\Ad}{Ad}
\DeclareMathOperator{\ad}{ad}
\DeclareMathOperator{\Exp}{Exp}

\begin{document}
\maketitle
\numberwithin{equation}{section}

\begin{abstract}
Let $\kkk$ be a local field, $G$ the set of $\kkk$-points of a connected semisimple algebraic $\kkk$-group $\mathbf{G}$, and $H$ the set of $\kkk$-points of a connected reductive algebraic $\kkk$-subgroup $\mathbf{H}$ of~$\mathbf{G}$ such that $\mathrm{rank}_{\kkk}(\mathbf{H})=\mathrm{rank}_{\kkk}(\mathbf{G})-1$.
We consider discrete subgroups $\Gamma$ of~$G$ acting properly discontinuously on $G/H$ and we examine their images under a Cartan projection $\mu : G\rightarrow V^+$, where $V^+$ is a closed convex cone in a real finite-dimensional vector space.
We show that if $\Gamma$ is neither a torsion group nor a virtually cyclic group, then $\mu(\Gamma)$ is almost entirely contained in one connected component of $V^+\setminus C_H$, where $C_H$ denotes the convex hull of $\mu(H)$ in~$V^+$.
As an application, we describe all torsion-free discrete subgroups of $G\times G$ acting properly discontinuously on $G$ by left and right translation when $\mathrm{rank}_{\kkk}(\mathbf{G})=1$.
\vspace{0.2cm}
\end{abstract}

\section{Introduction}

Let $\kkk$ be a local field, $G$ the set of $\kkk$-points of a connected semisimple algebraic $\kkk$-group of rank one, and $\Delta_G$ the diagonal of $G\times G$.
In this paper we describe all torsion-free discrete subgroups of $G\times G$ acting properly discontinuously on $(G\times G)/\Delta_G$ (Theorem~\ref{coro, rang un}).
To this end, we prove a general result on the Cartan projection of discrete groups acting properly discontinuously on corank-one reductive homogeneous spaces (Theorem~\ref{theoreme principal, general}).
This result holds for algebraic groups over any local field, but we first state it in the setting of real Lie groups (Theorem~\ref{theoreme principal, reel}).

\subsection{The main result in the real case}\label{Cas reel}
Let $G$ be a real connected semisimple linear Lie group and $H$ a closed connected reductive subgroup of~$G$.
It is known that $G$ contains an infinite discrete subgroup $\Gamma$ acting properly discontinuously on~$G/H$ if and only if $\mathrm{rank}_{\R}(H)<\mathrm{rank}_{\R}(G)$; this is the Calabi-Markus phenomenon (\cite{ko1}, Cor.~4.4). In this paper we consider the case when $\mathrm{rank}_{\R}(H)=\mathrm{rank}_{\R}(G)-1$.

Let us introduce some notation.
Fix a Cartan subgroup $A$ of~$G$ with Lie algebra~$\aaa$.
Denote by~$\Phi=\Phi(A,G)$ the system of restricted roots of~$A$ in~$G$, by~$\Phi^+$ a system of positive roots, by~$A^+ = \{ a\in A,\ \chi(a)\geq 1\ \, \forall\chi\in\Phi^+\} $ the corresponding closed positive Weyl chamber, and set $V^+ = \log A^+ \subset\nolinebreak \aaa$.
There is a maximal compact subgroup $K$ of~$G$ such that the Cartan decomposition $G=\nolinebreak KA^+K$ holds: every element $g\in G$ may be written as $g=k_1ak_2$ for some $k_1,k_2\in K$ and a unique $a\in A^+$ (\cite{hel}, Chap.~9, Th.~1.1). Setting $\mu(g)=\log a$ defines a map $\mu : G\rightarrow V^+$, which is continuous, proper, and surjective. It is called the \emph{Cartan projection} relative to the Cartan decomposition $G=KA^+K$.

Since $\mathrm{rank}_{\R}(H)=\mathrm{rank}_{\R}(G)-1$, the set $\mu(H)$ separates $V^+$ into finitely many connected components, which are permuted by the opposition involution~$\iota$.
(Recall that for every $a\in A^+$ we have $\iota(\log a)=\log a'$, where $a'$ is the unique element of~$A^+$ conjugate to~$a^{-1}$.)

In this setting our main result is the following.

\begin{theo}\label{theoreme principal, reel}
Let $G$ be a real connected semisimple linear Lie group and~$H$~a closed connected reductive subgroup of~$G$ such that $\mathrm{rank}_{\R}(H) = \mathrm{rank}_{\R}(G)-\nolinebreak 1$.
For~every discrete subgroup $\Gamma$ of~$G$ acting properly discontinuously on $G/H$, there exists a connected component $C$ of $V^+\setminus\mu(H)$ such that $\mu(\gamma)\in C\cup\iota(C)$ for almost all $\gamma\in\Gamma$.
If $\Gamma$ is not virtually cyclic, then $\iota(C)=C$.
\end{theo}

Recall that a group $\Gamma$ is said to satisfy some property \emph{virtually} if it contains a subgroup of finite index satisfying this property.
A property is said to be true for \emph{almost all} $\gamma\in\Gamma$ if it is true for all $\gamma\in\Gamma$ with at most finitely many exceptions.

By results of Chevalley (\cite{che}, Chap.~2, Th.~14 \& 15), if $G$ is a real connected semisimple linear Lie group and~$H$~a closed connected reductive subgroup of~$G$, then $G$ (resp.~$H$) is the identity component (for the real topology) of the set of $\R$-points of a connected semisimple linear algebraic $\R$-group~$\mathbf{G}$ (resp.~of a connected reductive algebraic $\R$-subgroup~$\mathbf{H}$ of~$\mathbf{G}$).
Theorem~\ref{theoreme principal, reel} is equivalent to the analogous result where $G$ (resp.~$H$) is replaced by~$\mathbf{G}(\R)$ (resp. by~$\mathbf{H}(\R)$).
We prove this result not only for $\R$-groups, but more generally for algebraic groups over any local field~$\kkk$.

\subsection{The main result in the general case}\label{Cas general}
Let $\kkk$ be a local field, \emph{i.e.}, $\R$, $\C$, a finite extension of~$\Q_p$, or the field $\F_q((t))$ of formal Laurent series over a finite field $\F_q$.
Let $G$ be the set of $\kkk$-points of a connected semisimple algebraic $\kkk$-group~$\mathbf{G}$ and $H$ the set of $\kkk$-points of a connected reductive algebraic $\kkk$-subgroup $\mathbf{H}$ of~$\mathbf{G}$ such that $\mathrm{rank}_{\kkk}(\mathbf{H}) = \mathrm{rank}_{\kkk}(\mathbf{G})-1$.
There is a Cartan projection~$\mu$ of~$G$ to a closed convex cone $V^+$ in some real finite-dimensional vector space (see Section~\ref{Decomposition de Cartan}).
The convex hull $C_H$ of~$\mu(H)$ in~$V^+$ separates~$V^+$ into finitely many connected components.
The opposition involution $\mu(G)\rightarrow\mu(G)$, which maps $\mu(g)$ to~$\mu(g^{-1})$ for all $g\in G$, extends to an involution~$\iota$ of~$V^+$ preserving $C_H$ and permuting the connected components of $V^+\setminus C_H$ (see Section~3.1).
Our main result in this general setting is the following.

\begin{theo}\label{theoreme principal, general}
Let $\kkk$ be a local field, $G$ the set of $\kkk$-points of a connected semisimple algebraic $\kkk$-group $\mathbf{G}$, and $H$ the set of $\kkk$-points of a connected reductive algebraic $\kkk$-subgroup~$\mathbf{H}$ of~$\mathbf{G}$ such that $\mathrm{rank}_{\kkk}(\mathbf{H}) = \mathrm{rank}_{\kkk}(\mathbf{G})-1$.
For~every discrete subgroup $\Gamma$ of~$G$ that acts properly discontinuously on~$G/H$ and that is not a torsion group, there exists a connected component~$C$ of $V^+\setminus C_H$ such that $\mu(\gamma)\in C\cup\iota(C)$ for almost all $\gamma\in\Gamma$.
If $\Gamma$ is not virtually cyclic, then $\iota(C)=C$.
\end{theo}

When $\kkk$ has characteristic zero, Theorem~\ref{theoreme principal, general} holds without assuming that $\Gamma$ is not a torsion group: indeed, in this case every discrete torsion subgroup of~$G$ is finite (Lemma~\ref{Sous-groupes discrets de torsion}).
This is not true when $\kkk=\F_q((t))$ for some finite field $\F_q$: in positive characteristic there are infinite discrete torsion subgroups of~$G$ that do not satisfy the conclusions of Theorem~\ref{theoreme principal, general}. We will give an example of such a group in Section~5.2.

\subsection{An application to $(G\times G)/\Delta_G$}
Our first application of Theorem~\ref{theoreme principal, general}, which is actually the main motivation of this paper, concerns homogeneous spaces of the form $(G\times G)/\Delta_G$, where $G$ is the set of $\kkk$-points of a connected semisimple algebraic $\kkk$-group $\mathbf{G}$ with $\mathrm{rank}_{\kkk}(\mathbf{G})=1$, and where $\Delta_G$ is the diagonal of $G\times G$.
In this situation, if $\mu$ is a Cartan projection of~$G$, then $\mu\times\mu$ is a Cartan projection of $G\times G$; we identify $V^+$ with $\R^+\times\R^+$ and $C_H$ with the diagonal of $\R^+\times\R^+$.

\begin{theo}\label{coro, rang un}
Let $\kkk$ be a local field, $G$ the set of $\kkk$-points of a connected semisimple algebraic $\kkk$-group $\mathbf{G}$ with $\mathrm{rank}_{\kkk}(\mathbf{G})=1$, and $\Delta_G$ the diagonal of~$G\times\nolinebreak G$.
Let $\Gamma$ be a discrete subgroup of $G\times G$.
\begin{enumerate}
	\item Assume that $\Gamma$ is torsion-free. Then it acts properly discontinuously on $(G\times G)/\Delta_G$ if and only if, up to switching the factors of $G\times G$, it is a graph of the form
$$\{ (\gamma,\varphi(\gamma)),\ \gamma\in\Gamma_0\} ,$$
where $\Gamma_0$ is a discrete subgroup of~$G$ and $\varphi : \Gamma_0\rightarrow G$ is a group homomorphism such that for all $R>0$, almost all $\gamma\in\Gamma_0$ satisfy $\mu(\varphi(\gamma))<\mu(\gamma)-R$.
  \item Assume that $\Gamma$ is residually finite and is not a torsion group. Then it acts properly discontinuously on $(G\times G)/\Delta_G$ if and only if, up to switching the factors of $G\times G$, it has a finite-index subgroup $\Gamma'$ that is a graph as in~\emph{(1)}.
\end{enumerate}
\end{theo}

Note that $(g,h)\Delta_G\mapsto gh^{-1}$ defines a $(G\times G)$-equivariant isomorphism from $(G\times G)/\Delta_G$ to~$G$, where $G\times G$ acts on~$G$ by $(g_1,g_2)\cdot g = g_1 g g_2^{-1}$.
Thus Theorem~\ref{coro, rang un} describes all torsion-free discrete subgroups of $G\times G$ acting properly discontinuously on~$G$ by left and right translation.

Recall that a group is said to be \emph{residually finite} if the intersection~of its normal finite-index subgroups is trivial.
It is known that if $\Gamma\subset G\times G$ is finitely generated, then it is residually finite (see~\cite{alp}, Cor.~1); if moreover $\kkk$ has characteristic zero, then $\Gamma$ has a finite-index subgroup that is torsion-free by Selberg's lemma (\cite{sel}, Lem.~8).

In the case of $G=\PSL_2(\R)$, Theorem~\ref{coro, rang un} has been proved for torsion-free groups by Kulkarni and Raymond~\cite{kr}.
In~\cite{ko2}, Kobayashi considered the more general case when $G$ is a real connected semisimple linear Lie group with $\mathrm{rank}_{\R}(G)=\nolinebreak 1$: he showed that every torsion-free discrete subgroup of $G\times G$ acting properly discontinuously on $(G\times G)/\Delta_G$ is a graph, and asked whether one of the two projections of this graph is always discrete in~$G$.
Theorem~\ref{coro, rang un} above answers this question positively and generalizes Kobayashi's result to all local fields.
It gives a complete description of all torsion-free discrete subgroups of $G\times\nolinebreak G$ acting properly discontinuously on $(G\times G)/\Delta_G$ in terms of a Cartan projection of~$G$.

Theorem~\ref{coro, rang un} applies to three-dimensional compact \emph{anti-de Sitter} manifolds, \emph{i.e.}, to three-dimensional compact Lorentz manifolds with constant sectional curvature~$-1$. Indeed, such manifolds are modeled on
$$\mathrm{AdS}^3 = \big\{ (x_1,x_2,x_3,x_4)\in\R^4,\quad x_1^2 + x_2^2 - x_3^2 - x_4^2 = 1\big\} $$
endowed with the Lorentz metric induced by $x_1^2+x_2^2-x_3^2-x_4^2$, which identifies with $(\SL_2(\R)\times\SL_2(\R))/\Delta_{\SL_2(\R)}$ (see Section~5.3).
Since three-dimensional compact anti-de Sitter manifolds are complete~\cite{kli}, they are quotients of the universal covering of $\mathrm{AdS}^3$.
By~\cite{kr}, up to a finite covering, they may in fact be written as
$$\Gamma\backslash (\PSL_2(\R)\times\PSL_2(\R))/\Delta_{\PSL_2(\R)},$$
where $\Gamma$ is a torsion-free discrete subgroup of $\PSL_2(\R)\times\PSL_2(\R)$ acting properly discontinuously on $(\PSL_2(\R)\times\PSL_2(\R))/\Delta_{\PSL_2(\R)}$.
We refer the reader to the introduction of~\cite{sal} for more details.

More generally, for any local field~$\kkk$ and any quadratic form $Q$ of Witt index two on~$\kkk^4$, the quadric
$$S(Q) = \{ x\in\kkk^4,\ Q(x)=1\} $$
identifies with $(\SL_2(\kkk)\times\SL_2(\kkk))/\Delta_{\SL_2(\kkk)}$ (see Section~5.3).
Theorem~\ref{coro, rang un} therefore applies to the discrete subgroups of $\SL_2(\kkk)\times\SL_2(\kkk)$ acting properly discontinuously on $S(Q)$.

Note that Theorem~\ref{coro, rang un} cannot be generalized to groups $G$ of higher rank. Indeed, take for instance $G = \SO(2,2n)$, and let $\Gamma_1$ (resp.\ $\Gamma_2$) be a torsion-free discrete subgroup of~$\SO(1,2n)$ (resp.\ of~$\U(1,n)$), where $\SO(1,2n)$ (resp.\ $\U(1,n)$) is seen as a subgroup of~$G$.
By \cite{ko1}, Prop.~4.9, $\Gamma_1\times\Gamma_2$ acts properly discontinuously on $(G\times G)/\Delta_G$. Other examples are obtained by replacing the triple $(\SO(2,2n),\SO(1,2n),\U(1,n))$ by $(\SO(4,4n),\SO(3,4n),\Sp(1,n))$ or by $(\U(2,2n),\U(1)\times\U(1,n),\Sp(1,n))$ (see~\cite{ko1}).

\subsection{An application to $\SLn(\kkk)/\SL_{n-1}(\kkk)$}
As another application of Theorem~\ref{theoreme principal, general}, we give a simpler proof of the following result, due to Benoist~\cite{be1}.

\begin{coro}\label{coro SLn}
Let $\kkk$ be a local field of characteristic zero.
If $n\geq 3$ is odd, then every discrete subgroup of~$\SLn(\kkk)$ acting properly discontinuously on $\SLn(\kkk)/\SL_{n-1}(\kkk)$ is virtually abelian.
\end{coro}

Theorem~\ref{theoreme principal, general} actually implies a slightly stronger version of Corollary~\ref{coro SLn}: we may replace ``virtually abelian'' by ``virtually cyclic''.

One consequence of Corollary~\ref{coro SLn} is that in characteristic zero if $n\geq\nolinebreak 3$ is odd, then the homogeneous space $\SLn(\kkk)/\SL_{n-1}(\kkk)$ has no compact quotient, \emph{i.e.}, there is no discrete subgroup $\Gamma$ of~$\SLn(\kkk)$ acting properly discontinuously on $\SLn(\kkk)/\SL_{n-1}(\kkk)$ with $\Gamma\backslash\SLn(\kkk)/\SL_{n-1}(\kkk)$ compact (see~\cite{be1}).

\subsection{Organization of the paper}
In Section~\ref{Decomposition de Cartan} we recall basic facts about Bruhat-Tits buildings, Cartan decompositions, and Cartan projections.
Section~\ref{Actions propres} is devoted to the proof of Theorem~\ref{theoreme principal, general}\,; we also discuss the assumption that $\Gamma$ is not a torsion group.
In Section~\ref{Application a SL_n/SL_{n-1}} we show how Theorem~\ref{theoreme principal, general} implies Corollary~\ref{coro SLn} in the case of~$G=\SLn(\kkk)$ and $H=\SL_{n-1}(\kkk)$.
In Section~\ref{(G*G)/Delta_G} we prove Theorem~\ref{coro, rang un}\,; we also show that the hypothesis that $\Gamma$ is not a torsion group is necessary in positive characteristic, and we describe our application to three-dimensional quadrics.

\subsection*{Acknowledgements}
I warmly thank Yves Benoist for his advice and encouragement.

\section{Cartan projections}\label{Decomposition de Cartan}

Throughout this article, we denote by~$\kkk$ a local field, \emph{i.e.}, $\R$, $\C$, a finite extension of~$\Q_p$, or the field $\F_q((t))$ of formal Laurent series over a finite field~$\F_q$.
If $\kkk=\R$ or $\C$, we denote by~$|\cdot|$ the usual absolute value on~$\kkk$; we set $\kkk^+=[1,+\infty[$.
If $\kkk$ is nonarchimedean, we denote by~$\mathcal{O}$ the ring of integers of~$\kkk$, by~$q$ the cardinal of the residue field of~$\kkk$, by~$\omega$ the (additive) valuation on~$\kkk$ sending any uniformizer to~$1$, and by $|\cdot| = q^{-\omega(\cdot)}$ the corresponding (multiplicative) absolute value; we set $\kkk^+=\{ x\in\kkk,\ |x|\geq 1\} $.
If $\mathbf{G}$ is an algebraic group, we denote by~$G$ the set of its $\kkk$-points and by~$\g$ the Lie algebra of~$G$.

In this section, we recall a few well-known facts on connected semisimple algebraic $\kkk$-groups and their Cartan projections.

\subsection{Weyl chambers}
Fix a connected semisimple algebraic $\kkk$-group $\mathbf{G}$.
Recall that the $\kkk$-split $\kkk$-tori of~$\mathbf{G}$ are all conjugate over~$\kkk$ (\cite{bot}, Th.~4.21).
Fix such a torus $\mathbf{A}$ and let $\mathbf{N}$ (resp.\ $\mathbf{Z}$) denote its normalizer (resp.\ centralizer) in~$\mathbf{G}$.
The group $X(\mathbf{A})$ of $\kkk$-characters of~$\mathbf{A}$ and the group $Y(\mathbf{A})$ of $\kkk$-cocharacters are both free $\Z$-modules of rank $r=\mathrm{rank}_{\kkk}(\mathbf{G})$, and there is a perfect pairing
$$\langle\cdot\,,\cdot\rangle : X(\mathbf{A})\times Y(\mathbf{A})\longrightarrow\Z.$$
If $\kkk$ is nonarchimedean, we set $A^{\circ}=A$; if $\kkk=\R$ or~$\C$, we set
$$A^{\circ} = \big\{ a\in A\,,\quad \chi(a)\in\, ]0,+\infty[ \quad\forall\chi\in X(\mathbf{A})\big\} .$$
The set $\Phi=\Phi(\mathbf{A},\mathbf{G})$ of restricted roots of~$\mathbf{A}$ in~$\mathbf{G}$, \emph{i.e.}, the set of nontrivial weights of~$\mathbf{A}$ in the adjoint representation of~$\mathbf{G}$, is a root system of the real vector space $V=Y(\mathbf{A})\otimes_{\Z}\R$ (\cite{bot}, Cor.~5.8).
The group $W=N/Z$ is finite and identifies with the Weyl group of~$\Phi$ (\cite{bot}, \S~5.1 \& Th.~5.3).
Choose a basis $\Delta = \{ \alpha_1,\ldots,\alpha_r\} $ of~$\Phi$ and let
$$\begin{array}{lcclcl}
& A^+ & = & \big\{ a\in A^{\circ}\!, & \ \alpha_i(a)\in\kkk^+ & \forall\,1\leq i\leq r \big\} \\
\mathrm{\big(resp.}\quad & V^+ & = & \big\{ x\in V, & \langle\alpha_i,x\rangle\geq 0 & \forall\,1\leq i\leq r \big\} \mathrm{\big)}
\end{array}$$
denote the closed positive Weyl chamber in~$A^{\circ}$ (resp.\ in~$V$) corresponding to~$\Delta$; the set $V^+$ is a closed convex cone in~$V$.
If $\kkk=\R$ or~$\C$, then $V$ identifies with~$\aaa$ and~$V^+$ with $\log A^+\subset\aaa$; we endow $V$ with the Euclidean norm $\Vert\cdot\Vert$ induced by the Killing form of $\g$.
If $\kkk$ is nonarchimedean, we endow $V$ with any $W$-invariant Euclidean norm $\Vert\cdot\Vert$.

\subsection{The Bruhat-Tits building}
In this subsection we assume $\kkk$ to be nonarchimedean. We briefly recall the construction of the Bruhat-Tits building of~$G$, which is a metric space on which $G$ acts properly discontinuously by isometries with a compact fundamental domain.
We refer to the original articles \cite{bt1} and \cite{bt2}, but the reader may also find \cite{rou} useful.

Let $\Res$ denote the restriction homomorphism from $X(\mathbf{Z})$ to $X(\mathbf{A})$, where $X(\mathbf{Z})$ denotes the group of $\kkk$-characters of~$\mathbf{Z}$.
There is a unique group homomorphism $\nu : Z\rightarrow V$ such that
$$\langle\Res(\chi),\nu(z)\rangle = -\,\omega(\chi(z))$$
for all $\chi\in X(\mathbf{Z})$ and $z\in Z$.
The set $\nu(Z)$ is a lattice in~$V$, and $\nu(A)$ is a sublattice of~$\nu(Z)$ of finite index.
The action of~$Z$ on~$V$ by translation along~$\nu(Z)$ extends to an action of~$N$ on~$V$ by affine isometries; such an extension is unique up to translation.

For every $\alpha\in\Phi$, let $\mathbf{U}_{\alpha}$ denote the connected unipotent $\kkk$-subgroup of~$\mathbf{G}$ corresponding to the root $\alpha$, as defined in~\cite{bt2}; the Lie algebra of~$U_{\alpha}$ is $\g_{\alpha}\oplus\g_{2\alpha}$, where $\g_{i\alpha}$ is the subspace of elements $X\in\g$ such that $\Ad(a)(X)=\alpha(a)^i X$ for all $a\in A$.
For every $u\in U_{\alpha}$, $u\neq 1$, the set $N\cap U_{-\alpha}\,u\,U_{-\alpha}$ has a unique element, which acts on~$V$ by the orthogonal reflection in some affine hyperplane~$\mathcal{H}_u$, defined by an equation of the form $\langle\alpha,x\rangle+\psi_{\alpha}(u)=0$, where $\psi_{\alpha}(u)\in\R$.
For every $x\in V$, set
$$U_{\alpha,x} = \big\{ u\in U_{\alpha},\quad u=1\ \ \mathrm{or}\ \ \langle\alpha,x\rangle+\psi_{\alpha}(u)\geq 0\big\} ;$$
by~\cite{bt2} it is a subgroup of~$U_{\alpha}$.
Set $N_x = \{ n\in N,\ n\cdot x=x\} $ and let $K_x$ denote the subgroup of~$G$ generated by $N_x$ and the subgroups $U_{\alpha,x}$, where $\alpha\in\Phi$.
The group $K_x$ is a maximal compact open subgroup of~$G$.

With this notation, the \emph{Bruhat-Tits building} $X$ of~$G$ is the set of equivalence classes of $G\times V$ for the relation
$$(g,x)\sim (g',x') \quad\Longleftrightarrow\quad \exists\,n\in N\ \mathrm{such\ that}\ x'=n\cdot x\ \mathrm{et}\ g^{-1}g'n\in K_x.$$
We endow $X$ with the quotient topology induced by the discrete topology of~$G$ and the Euclidean structure of~$V$.
By construction, $V$ embeds into~$X$; we identify~it with its image in~$X$.
The group $G$ acts on~$X$ by
$$g'\cdot\overline{(g,x)}\ =\ \overline{(g'g\,,x)},$$
where $\overline{(g,x)}$ denotes the image of $(g,x)\in G\times V$ in~$X$.
This action is properly discontinuous, with a compact fundamental domain.
By construction, the stabilizer of any point~$x\in V$ is~$K_x$.
The \emph{apartments} of~$X$ are the sets $g\cdot V$, where $g\in G$; the \emph{walls} of~$X$ are the sets $g\cdot\mathcal{H}_u$, where $g\in G$ and $u\in U_{\alpha}$ for some $\alpha\in\Phi$.
A \emph{chamber} of~$X$ (or~\emph{alcove}) is a connected component of~$X$ deprived of its walls.
The space $X$ has the following property: for any pair $(x,x')$ of points in~$X$, there is an apartment containing both $x$ and~$x'$. We can therefore endow~$X$ with a distance~$d$ defined as follows: $d(x,x')$ is the Euclidean distance between $x$ and~$x'$ in any apartment containing $x$ and $x'$ (it does not depend on the apartment). The group~$G$ acts on~$X$ by isometries for this distance.

\subsection{Cartan decompositions and Cartan projections}\label{Projection de Cartan}
If $\kkk=\R$ or $\C$, then there is a maximal compact subgroup $K$ of~$G$ such that the Cartan decomposition $G=KA^+K$ holds: for every $g\in G$, there are elements $k_1,k_2\in K$ and a unique $a\in A^+$ such that $g = k_1 a k_2$ (\cite{hel}, Chap.~9, Th.~1.1).
Setting $\mu(g)=\log a$ defines a map $\mu : G\rightarrow V^+\simeq\log A^+$, which is continuous, proper, and surjective.
It is called the \emph{Cartan projection} relative to the Cartan decomposition $G=KA^+K$.

Now assume $\kkk$ to be nonarchimedean.
Consider the extremal point $x_0$ of the closed cone~$V^+$, defined by $\langle\alpha_i,x_0\rangle=0$ for all $1\leq i\leq r$, and set $K=K_{x_0}$.
Let $Z^+\subset Z$ denote the inverse image of~$V^+$ under~$\nu$.
By~\cite{bt1} the group $G$ acts transitively on the set of couples $(\mathcal{A},\mathcal{C})$, where $\mathcal{A}$ is an apartment of~$X$ and $\mathcal{C}$~is a chamber of~$X$ contained in~$\mathcal{A}$.
This can be translated into algebraic terms by the existence of a \emph{Cartan decomposition} $G = KZ^+K$: for every $g\in G$ there are elements $k_1$, $k_2\in K$ and $z\in Z^+$ such that $g = k_1 z k_2$, and $\nu(z)$ is uniquely defined.
Setting $\mu(g)=\nu(z)$ defines a map $\mu : G\rightarrow V^+$, which is continuous and proper; its image $\mu(G)$ is the intersection of $V^+$ with a lattice of $V$.
The map $\mu$ is called the \emph{Cartan projection} relative to the Cartan decomposition $G=KZ^+K$.

\subsection{A geometric interpretation}\label{Interpretation geometrique}
Let $X$ be either the Riemannian symmetric space $G/K$ if $\kkk=\R$ or $\C$, or the Bruhat-Tits building of~$G$ if $\kkk$ is nonarchimedean.
We now recall a geometric interpretation of the Cartan projection~$\mu$ in terms of a distance on~$X$.

Assume that $\kkk=\R$ or $\C$, and let $\g = \kk + \p$ be the Cartan decomposition of~$\g$ corresponding to the Cartan decomposition $G = KA^+K$.
The Killing form~$\kappa$ of~$\g$ is definite positive on~$\p$, hence induces a Euclidean norm $\Vert\cdot\Vert$ on $\p$.
Let $\pi$ denote the natural projection of $G$ onto $X=G/K$, and set $x_0=\pi(1)\in X$.
The map~$d\pi_1$ realizes an isomorphism between~$\p$ and the tangent space of~$X$ at~$x_0$; thus $\kappa|_{\p\times\p}$ induces a $G$-invariant Riemannian metric on~$X$.
Let~$d$ denote the corresponding distance on~$X$.
The following result is probably well known; we prove it for the reader's convenience.

\begin{lem}[$\kkk=\R$ or $\C$]\label{rho contraction, reel}
Let $\rho : X\rightarrow V^+$ denote the map sending $x=g\cdot x_0\in X$ to~$\mu(g)$.
For all $x,x'\in X$,
$$\Vert\rho(x)-\rho(x')\Vert \leq d(x,x').$$
Moreover, the restriction of $\rho$ to $A^+\cdot x_0$ is an isometry.
\end{lem}

\begin{proof}
We identify $V^+$ with $\log A^+\subset\aaa$.
Let $\Exp : \p\rightarrow X$ denote the exponential diffeomorphism mapping $Y\in\p$ to~$\gamma_Y(1)$, where $\gamma_Y$ is the unique geodesic in~$X$ such that $\gamma_Y(0)=\nolinebreak x_0$ and $\gamma'_Y(0)=d\pi_1(Y)$.
For every $x\in X$, there exists $k\in K$ such that $x=k\exp(\rho(x))\cdot x_0$; by \cite{hel}, Chap.~4, Th.~3.3,
\begin{equation}\label{Exp et rho}
x = \Exp\big((\Ad k)(\rho(x))\big).
\end{equation}
Fix $x,x'\in X$ and let $\gamma=(y_t)_{t\in [0,1]}$ be the geodesic segment from $y_0=\nolinebreak x$ to $y_1=\nolinebreak x'$.
By \cite{hel}, p.~295, and~(\ref{Exp et rho}), the map $t\mapsto\rho(y_t)$ is smooth and there exists a smooth map $t\mapsto k_t$ from $[0,1]$ to~$K$ such that $y_t = \Exp((\Ad k_t)(\rho(y_t)))$ for all $t\in [0,1]$.
Since~$X$ has nonpositive sectional curvature (\cite{hel}, Chap.~5, Th.~3.1), the length of~$\gamma$ in $X$ is not less than the length of~$\Exp^{-1}(\gamma)$ in~$\p$ (\cite{hel}, Chap.~1, Th.~13.1), namely,
\begin{equation}\label{inegalite courbure}
d(x,x') \geq \int_0^1 \bigg\Vert\frac{\mathrm{d}\big((\Ad k_t)(\rho(y_t))\big)}{\mathrm{d}t}(t')\bigg\Vert\ \mathrm{d}t'.
\end{equation}
Now for all $t'\in [0,1]$,
$$\frac{\mathrm{d}\big((\Ad k_t)(\rho(y_t))\big)}{\mathrm{d}t}(t') = (\Ad k_{t'})\Big(\frac{\mathrm{d}(\rho(y_t))}{\mathrm{d}t}(t')\Big) + \Big(\frac{\mathrm{d}(\Ad k_t)}{\mathrm{d}t}(t')\Big)\big(\rho(y_{t'})\big),$$
where
$$(\Ad k_{t'})\Big(\frac{\mathrm{d}(\rho(y_t))}{\mathrm{d}t}(t')\Big) \in (\Ad k_{t'})(\aaa)$$
and
\begin{eqnarray*}
\Big(\frac{\mathrm{d}(\Ad k_t)}{\mathrm{d}t}(t')\Big)\big(\rho(y_{t'})\big) & = & (\Ad k_{t'})\bigg(\ad\Big(\frac{\mathrm{d}(k_{t'}^{-1}k_{t'+t})}{\mathrm{d}t}(0)\Big)\big(\rho(y_{t'})\big)\bigg)\\
& \in & (\Ad k_{t'})\big([\kk,\aaa]\big).
\end{eqnarray*}
The subspaces $\aaa$ and $[\kk,\aaa]$ are orthogonal with respect to~$\kappa$.
Indeed, the decomposition of~$\g$ into eigenspaces under the adjoint action of~$\aaa$ is orthogonal with respect to~$\kappa$ (\cite{hel}, Chap.~3, Th.~4.2); in particular, $\aaa$ is orthogonal to the sum~$[\g,\aaa]$ of the root spaces of~$\g$.
Since $\kappa$ is invariant under~$\Ad G$ (\cite{hel}, p.~131), the subspaces $(\Ad k_{t'})(\aaa)$ and $(\Ad k_{t'})([\kk,\aaa])$ are orthogonal with respect to~$\kappa$ and
\begin{eqnarray}
\bigg\Vert\frac{\mathrm{d}\big((\Ad k_t)(\rho(y_t))\big)}{\mathrm{d}t}(t')\bigg\Vert & \geq & \Big\Vert(\Ad k_{t'})\Big(\frac{\mathrm{d}(\rho(y_t))}{\mathrm{d}t}(t')\Big)\Big\Vert\label{derivee}\\
& = & \Big\Vert\frac{\mathrm{d}(\rho(y_t))}{\mathrm{d}t}(t')\Big\Vert.\nonumber
\end{eqnarray}
Thus
$$d(x,x')\ \geq\ \int_0^1 \Big\Vert\frac{\mathrm{d}(\rho(y_t))}{\mathrm{d}t}(t')\Big\Vert\ \mathrm{d}t'\ =\ \Vert\rho(x)-\rho(x')\Vert.$$
If $x,x'\in A^+\cdot x_0$, then $k_t=1$ for all $t\in [0,1]$; hence (\ref{derivee}) is an equality.
Moreover, in this case (\ref{inegalite courbure}) is also an equality since the geodesic submanifold $A\cdot x_0=\Exp(\aaa)$ has zero sectional curvature (\cite{hel}, Chap.~5, \S~3, Rem.~2).
This implies $d(x,x') = \Vert\rho(x)-\nolinebreak\rho(x')\Vert$.
\end{proof}

Since~$K$ fixes $x_0$ and since~$G$ acts on~$X$ by isometries, Lemma~\ref{rho contraction, reel} implies that for every $a\in A^+$ and every $g\in KaK$,
\begin{equation}\label{mu distance, reel}
d(g\cdot x_0,x_0) = d(a\cdot x_0,x_0) = \Vert\rho(a\cdot x_0)-\rho(x_0)\Vert = \Vert\mu(g)\Vert.
\end{equation}

Now assume $\kkk$ to be nonarchimedean and let $X$ denote the Bruhat-Tits building of~$G$, endowed with the distance~$d$ defined in Section~2.2.
Recall that $K=K_{x_0}$ is the stabilizer of the point $x_0\in V$ defined by $\langle\alpha_i,x_0\rangle=0$ for all $1\leq i\leq r$.
Since~$G$ acts on~$X$ by isometries and since~$V$ is isometrically embedded as an apartment in~$X$, for every $z\in Z^+$ and every $g\in KzK$,
\begin{equation}\label{mu distance, ultrametrique}
d(g\cdot x_0,x_0) = d(z\cdot x_0,x_0) = d(\mu(g),x_0) = \Vert\mu(g)\Vert,
\end{equation}
where $\Vert\cdot\Vert$ is the Euclidean norm on~$V$.
Lemma~\ref{rho contraction, reel} also holds in this setting.

\begin{lem}[$\kkk$ nonarchimedean]\label{rho contraction, ultrametrique}
Let $\rho : X\rightarrow V^+$ denote the map sending $x=g\cdot x_0\in X$ to~$\mu(g)$.
For all $x,x'\in X$,
$$\Vert\rho(x)-\rho(x')\Vert \leq d(x,x').$$
\end{lem}

\begin{proof}
Let~$\mathcal{C}$ denote the unique chamber in~$V^+$ containing~$x_0$.
We first recall the construction of a retraction $\rho_{V,\mathcal{C}} : X\rightarrow V$, as defined in \cite{bt1}, \S~2.3. \linebreak
For every $x\in X$, there is an apartment $\mathcal{A}$ containing both $x$ and~$\mathcal{C}$ (\cite{bt1}, Prop.~2.3.1), and there is an element $k\in K$ fixing $\mathcal{C}$ pointwise and mapping~$\mathcal{A}$ to~$V$ (\cite{bt1}, Prop.~2.3.2).
The point $k\cdot x\in V$ does not depend on the choice of~$\mathcal{A}$ and~$k$.
Setting $\rho_{V,\mathcal{C}}(x)=k\cdot x$ defines a map $\rho_{V,\mathcal{C}} : X\rightarrow V$ such that for all~$x,x'\in X$,
$$\Vert\rho_{V,\mathcal{C}}(x)-\rho_{V,\mathcal{C}}(x')\Vert\ \leq\ d(x,x')$$
(\cite{bt1}, Prop.~2.5.3).
We claim that for all~$x,x'\in X$,
\begin{equation}\label{les deux contractions}
\Vert\rho(x)-\rho(x')\Vert\ \leq\ \Vert\rho_{V,\mathcal{C}}(x)-\rho_{V,\mathcal{C}}(x')\Vert.
\end{equation}
Indeed, it follows from the definitions of $\rho$ and $\rho_{V,\mathcal{C}}$ that $\rho_{V,\mathcal{C}}(x)\in W\cdot\rho(x)$ for all $x\in X$.
Since the norm $\Vert\cdot\Vert$ is $W$-invariant, it is enough to show that
\begin{equation}\label{inegalite avec w}
\Vert\rho(x)-\rho(x')\Vert\ \leq\ \Vert \rho(x)-w\cdot\rho(x')\Vert
\end{equation}
for all $x,x'\in X$ and all $w\in W$.
Recall that $W$ is generated by the set~$S$ of orthogonal reflections in the hyperplanes $\{ x\in V,\ \langle\alpha_i,x\rangle = 0\big\} $, where $1\leq i\leq r$.
Write $w=s_m\ldots s_1$, where $s_j\in S$ for all~$j$.
We argue by induction on $m$.
If $(s_m\ldots s_1)\cdot\rho(x')\in V^+$, then $s_m\ldots s_1=1$ and (\ref{inegalite avec w}) is obvious.
Otherwise, the points $\rho(x)$ and $(s_m\ldots s_1)\cdot\nolinebreak\rho(x')$ lie in two distinct connected components of $V\setminus\mathcal{H}$, where $\mathcal{H}$ denotes the hyperplane of fixed points of~$s_m$.
Let $y$ be the intersection point of $\mathcal{H}$ with the line segment $[\rho(x),(s_m\ldots s_1)\cdot\rho(x')]$.
Since $s_m$ is an orthogonal reflection,
\begin{eqnarray*}
\Vert\rho(x)-(s_m\ldots s_1)\cdot\rho(x')\Vert & = & \Vert\rho(x)-y\Vert + \Vert y-(s_m\ldots s_1)\cdot\rho(x')\Vert\\
& = & \Vert\rho(x)-y\Vert + \Vert y-(s_{m-1}\ldots s_1)\cdot\rho(x')\Vert\\
& \geq & \Vert\rho(x)-(s_{m-1}\ldots s_1)\cdot\rho(x')\Vert.
\end{eqnarray*}
By the induction assumption, $\Vert\rho(x)-(s_m\ldots s_1)\cdot\rho(x')\Vert \geq \Vert\rho(x)-\rho(x')\Vert$.
This proves (\ref{les deux contractions}) and completes the proof of Lemma~\ref{rho contraction, ultrametrique}.
\end{proof}

The following result will be needed in the proof of Theorem~\ref{theoreme principal, general}.

\begin{lem}\label{inegalite pour mu}
Let $\kkk$ be a local field, $G$ the set of $\kkk$-points of a connected semisimple algebraic $\kkk$-group, and $\mu : G\rightarrow V^+$ a Cartan projection. For all $g,g'\in G$, the following two inequalities hold:
\begin{eqnarray}
\Vert\mu(gg')-\mu(g)\Vert & \leq & \Vert\mu(g')\Vert,\label{inegalite 1}\\
\Vert\mu(gg')-\mu(g')\Vert & \leq & \Vert\mu(g)\Vert.\label{inegalite 2}
\end{eqnarray}
\end{lem}

\begin{proof}
Since $G$ acts on~$X$ by isometries, (\ref{inegalite 1}) follows immediately from Lemmas~\ref{rho contraction, reel} and~\ref{rho contraction, ultrametrique}, together with Formulas~(\ref{mu distance, reel}) and~(\ref{mu distance, ultrametrique}).
We claim that~(\ref{inegalite 1}) implies~(\ref{inegalite 2}).
Indeed, if $w\in W$ denotes the ``longest'' element of~$W$, such that $w\cdot z^{-1}\in Z^+$ for all $z\in Z^+$, then $\mu(g^{-1})=w\cdot(-\mu(g))$ for all $g\in G$.
Since the norm~$\Vert\cdot\Vert$ on~$V$ is $W$-invariant, the opposition involution $\iota : \mu(G)\rightarrow\mu(G)$, which maps $\mu(g)$ to $\mu(g^{-1})$ for all $g\in G$, is an isometry.
Together with~(\ref{inegalite 1}), this implies
$$\Vert\mu(gg')-\mu(g')\Vert\ =\ \Vert\mu({g'}^{-1}g^{-1})-\mu({g'}^{-1})\Vert\ \leq\ \Vert\mu(g^{-1})\Vert\ =\ \Vert\mu(g)\Vert. \qedhere$$
\end{proof}

\section{Proper actions on $G/H$ in the corank-one case}\label{Actions propres}

In this section we give a proof of Theorem~\ref{theoreme principal, general} and we discuss the assumption that $\Gamma$ is not a torsion group.

\subsection{Proof of Theorem~\ref{theoreme principal, general}}\label{Demonstration du theoreme principal}
With the notation of Section~\ref{Decomposition de Cartan}, let $\mathbf{H}$ be a connected reductive algebraic $\kkk$-subgroup of~$\mathbf{G}$ with $\mathrm{rank}_{\kkk}(\mathbf{H})=\mathrm{rank}_{\kkk}(\mathbf{G})-\nolinebreak 1$.
Fix a maximal $\kkk$-split $\kkk$-torus~$\mathbf{A}_{\mathbf{H}}$ of~$\mathbf{H}$.
After conjugating~$\mathbf{H}$ by an element of~$G$, we may assume that $\mathbf{A}_{\mathbf{H}}\subset\mathbf{A}$ (\cite{bot}, Th.~4.21).
Recall that $\mathbf{H}$ is the almost product of a central torus and of its derived group, which is semisimple (\cite{bot}, Prop.~2.2). Therefore $H$ admits a Cartan decomposition $H = K_H Z_H^+ K_H$, where $\mathbf{Z}_{\mathbf{H}}$ is the centralizer of~$\mathbf{A}_{\mathbf{H}}$ in~$\mathbf{H}$ and $K_H$ is some maximal compact subgroup of~$H$.
We now use a result proved by Mostow~\cite{mos} and Karpelevich~\cite{kar} in the real case, and by Landvogt~\cite{lan} in the nonarchimedean case: after conjugating $\mathbf{H}$ by an element~of~$G$, we may assume that $K_H\subset\nolinebreak K$.
Thus $\mu(H)=\mu(Z_H)$ and the convex hull~$C_H$ of~$\mu(H)$ in~$V^+$ is the intersection~of $V^+$ with a finite union of hyperplanes of~$V$ parametrized by the Weyl group~$W$.
The opposition involution $\iota : \mu(G)\rightarrow\mu(G)$, which maps $\mu(g)$ to $\mu(g^{-1})$ for all $g\in G$, extends to an isometry of~$V^+$, still denoted by~$\iota$. It preserves $\mu(H)$, hence $C_H$, and permutes the connected components of~$V^+\setminus C_H$.

Our proof of Theorem~\ref{theoreme principal, general} is based on the \emph{properness criterion} of Benoist (\cite{be1}, Cor.~5.2) and Kobayashi (\cite{ko3}, Th.~1.1), which states that a subgroup $\Gamma$ of~$G$ acts properly discontinuously on~$G/H$ if and only if the set $\mu(\Gamma)\cap(\mu(H)+C')$ is bounded for every compact subset $C'$ of~$V$. This condition is equivalent to the boundedness of $\mu(\Gamma)\cap(C_H+C')$ for every compact subset~$C'$ of~$V$.

Our proof is also based on the following observation $(\ast)$: if $(x_n)_{n\in\N}$ is a sequence of points of~$V^+$ whose distance to~$C_H$ is larger than a given $R>0$, and if $\Vert x_{n+1}-x_n\Vert\leq R$ for all $n\in\N$, then all elements $x_n$ belong to the same connected component of~$V^+\setminus C_H$.

We now give a proof of Theorem~\ref{theoreme principal, general}.
Let $C_1,\ldots,C_s$ be the connected components of $V^+\setminus C_H$ and let $\Gamma$ be a discrete subgroup of~$G$ acting properly discontinuously on~$G/H$.
The set $\mu(\Gamma)$ is invariant under the opposition involution~$\iota$.

Assume that $\Gamma$ is not a torsion group and fix an element $\gamma\in\Gamma$ of infinite order.
Since~$\Gamma$ is discrete and since $\mu$ is a proper map, the sequence $(\Vert\mu(\gamma^n)\Vert)_{n\in\Z}$ tends to infinity as $n$ tends to~$\pm\infty$.
Let $F$ be the set of elements $\gamma'\in\Gamma$ such that the distance of~$\mu(\gamma')$ to~$C_H$ is $\leq\Vert\mu(\gamma)\Vert$.
From the discreteness of~$\Gamma$, the properness of~$\mu$, and the properness criterion, we deduce that $F$ is finite.
Moreover, by Lemma~\ref{inegalite pour mu},
$$\big\Vert\mu(\gamma^{n+1}) - \mu(\gamma^n)\big\Vert \leq \Vert\mu(\gamma)\Vert$$
for all $n\in\Z$.
By the observation $(\ast)$ above, there are integers $1\leq i,j\leq s$ such that $\mu(\gamma^n)\in C_i$ (resp.\ $\mu(\gamma^{-n})\in C_j$) for almost all $n\in\N$.
The opposition involution~$\iota$ interchanges $C_i$ and~$C_j$.

Note that for every $\gamma'\in\Gamma$, Lemma~\ref{inegalite pour mu} implies
$$\big\Vert\mu(\gamma'\gamma^n) - \mu(\gamma^n)\big\Vert \leq \Vert\mu(\gamma')\Vert$$
for all $n\in\Z$. By the properness criterion, $\mu(\gamma'\gamma^n)\in C_i$ and $\mu(\gamma'\gamma^{-n})\in C_j$ for almost all $n\in\N$.

First consider the case $i=j$. Let $F'$ be the set of elements $\gamma'\in\Gamma$ such that $\mu(\gamma')\notin C_i$.
We claim that $F'$ is finite.
Indeed, let $\gamma'\in F'$.
By Lemma~\ref{inegalite pour mu},
$$\big\Vert\mu(\gamma'\gamma^{n+1}) - \mu(\gamma'\gamma^n)\big\Vert \leq \Vert\mu(\gamma)\Vert$$
for all $n\in\Z$.
Moreover, $\mu(\gamma')\notin C_i$, and we have just seen that $\mu(\gamma'\gamma^n)\in C_i$ for almost all $n\in\Z$.
By the observation $(\ast)$ above, there is an integer $n\in\Z$ such that $\gamma'\gamma^n\in F$.
Therefore, $F'\subset F\gamma^{\Z}$.
Since~$F$ is finite and since for every $f\in F$ the element $f\gamma^n$ belongs to~$C_i$ for almost all $n\in\Z$, the set $F'$ is finite.
This proves the claim.

Now consider the case $i\neq j$.
We claim that the subgroup $\gamma^{\Z}$ has finite index in~$\Gamma$.
Indeed, let $\gamma'\in\Gamma$.
By Lemma~\ref{inegalite pour mu},
$$\big\Vert\mu(\gamma'\gamma^{n+1}) - \mu(\gamma'\gamma^n)\big\Vert \leq \Vert\mu(\gamma)\Vert$$
for all $n\in\Z$.
Moreover, we have seen that $\mu(\gamma'\gamma^n)\in C_i$ and $\mu(\gamma'\gamma^{-n})\in C_j$ for almost all $n\in\N$. By the observation $(\ast)$ above, there is an integer $n\in\Z$ such that $\gamma'\gamma^n\in F$.
Therefore, $\Gamma = F \gamma^{\Z}$.
Since~$F$ is finite, $\gamma^{\Z}$ has finite index in~$\Gamma$.
This proves the claim and completes the proof of Theorem~\ref{theoreme principal, general}.

\subsection{Discrete torsion groups in characteristic zero}\label{Torsion}
In this subsection we show that when $\kkk$ has characteristic zero, the assumption that $\Gamma$ is not a torsion group may be removed from Theorem~\ref{theoreme principal, general}.
When $\Gamma$ is known to be finitely generated, this follows from Selberg's lemma (\cite{sel}, Lem.~8).
In general it is also true, based on the following lemma, which is probably well~known.

\begin{lem}\label{Sous-groupes discrets de torsion}
Let $\kkk$ be a local field of characteristic zero and $\mathbf{G}$ a linear algebraic $\kkk$-group.
If $\kkk$ is a $p$-adic field, then every torsion subgroup of~$G$ is finite.
If $\kkk=\R$ or $\C$, then every discrete torsion subgroup of~$G$ is finite.
\end{lem}

\begin{proof}
Embed $\mathbf{G}$ in~$\mathbf{GL}_n$ for some $n\geq 1$.
Let $\Gamma$ be a torsion subgroup of~$G$.
By a result of Schur (\cite{cr}, Th.~36.14), $\Gamma$ contains a finite-index abelian subgroup $\Gamma'$ whose elements are all semisimple.
To show that $\Gamma$ is finite, it is enough to prove the finiteness of~$\Gamma'$.

Assume that $\kkk$ is a $p$-adic field.
The elements of~$\Gamma'$ are diagonalizable in a common basis over an algebraic closure of~$\kkk$.
For every $\gamma\in\Gamma'$ the eigenvalues of~$\gamma$ are roots of unity; they generate a cyclotomic extension~$\kkk_{\gamma}$ of~$\kkk$, and $[\kkk_{\gamma}:\kkk]\leq n$ since the characteristic polynomial of~$\gamma$ has degree~$n$.
Now there are only finitely many cyclotomic extensions of $\kkk$ of degree~$\leq n$ (\cite{neu}, Chap.~2, Th.~7.12 \& Prop.~7.13).
Therefore the field generated by all extensions $\kkk_{\gamma}$, $\gamma\in\Gamma'$, has finite degree over~$\kkk$, hence contains only finitely many roots of unity (\cite{neu}, Chap.~2, Prop.~5.7).
This implies the finiteness of~$\Gamma'$.

Assume that $\kkk=\R$ or $\C$ and that in addition $\Gamma$ is discrete in~$G$.
The elements of~$\Gamma'$ are diagonalizable in a common basis over~$\C$, and their eigenvalues are roots of unity.
Since the group $\mathbb{U}$ of complex numbers of modulus one is compact, every discrete subgroup of~$\mathbb{U}^n$ is finite.
This implies the finiteness of~$\Gamma'$.
\end{proof}

When $\kkk$ has positive characteristic, there exist infinite discrete torsion subgroups in~$G$.
They all have a unipotent subgroup of finite index (this follows from \cite{tits}, Prop.~2.8, for instance).
Some of them do not satisfy the conclusions of Theorem~\ref{theoreme principal, general}\,: we will give an example of such a group in Section~5.2.

\section{An application to $\SLn(\kkk)/\SL_{n-1}(\kkk)$}\label{Application a SL_n/SL_{n-1}}

In this section we discuss the case of $G=\SLn(\kkk)$ and $H=\SL_{n-1}(\kkk)$.
We show how Theorem~\ref{theoreme principal, general} implies Corollary~\ref{coro SLn}.

Let $\mathbf{G}=\mathbf{SL}_n$ for some integer $n\geq 2$. The group $\mathbf{A}$ of diagonal matrices of determinant one is a maximal $\kkk$-split $\kkk$-torus of~$\mathbf{G}$, which is its own centralizer, \emph{i.e.}, $\mathbf{Z}=\mathbf{A}$.
The corresponding roots are the linear forms $\varepsilon_i-\varepsilon_j$, $1\leq i\neq j\leq n$, where
$$\varepsilon_i\big(\mathrm{diag}(a_1,\ldots,a_n)\big) = a_i.$$
A basis of the root system of~$\mathbf{A}$ in~$\mathbf{G}$ is given by the roots $\varepsilon_i-\varepsilon_{i+1}$, where $1\leq i\leq n-1$.
If $\kkk=\R$ or~$\C$ (resp. if $\kkk$ is nonarchimedean), the corresponding positive Weyl chamber is
\begin{eqnarray*}
A^+ & = & \big\{ \mathrm{diag}(a_1,\ldots,a_n)\in A,\ \ \! a_i\in\, ]0,+\infty[\ \forall i\ \,\mathrm{and}\ a_1\geq\ldots\geq a_n\big\} \\
\mathrm{\big(resp.}\quad A^+ & = & \big\{ \mathrm{diag}(a_1,\ldots,a_n)\in A,\ |a_1|\geq\ldots\geq|a_n|\big\} \mathrm{\big).}
\end{eqnarray*}
Set $K=\SO(n)$ (resp. $K=\SU(n)$, resp. $K=\SLn(\mathcal{O})$) if $\kkk=\R$ (resp. if \linebreak $\kkk=\C$, resp. if $\kkk$ is nonarchimedean).
The Cartan decomposition $G=\nolinebreak KA^+K$ holds.
If $\kkk=\R$ (resp.\ if $\kkk=\nolinebreak\C$), it follows from the polar decomposition in $\GLnR$ (resp.\ in $\GLnC$) and from the reduction of symmetric (resp.\ of Hermitian) matrices; if~$\kkk$ is nonarchimedean, it follows from the structure theorem for finitely generated modules over a principal ideal domain.
The real vector space
$$V = \big\{ (x_1,\ldots,x_n)\in\R^n,\ x_1+\ldots+x_n=0\big\}\ \simeq\ \R^{n-1}$$
and its closed convex cone
$$V^+ = \big\{ (x_1,\ldots,x_n)\in V,\ x_1\geq\ldots\geq x_n\big\} $$
do not depend on~$\kkk$.
Let $\mu : G\rightarrow V^+$ denote the Cartan projection relative to the Cartan decomposition~$G=KA^+K$.
If $\kkk=\R$ or~$\C$, then $\mu(g)=(x_1,\ldots,x_n)$, where $e^{2x_i}$ is the $i$-th eigenvalue of~$^t\!\overline{g}g$.

Let $\mathbf{H}=\mathbf{SL}_{n-1}$, which we consider as a subgroup of~$\mathbf{G}$ by embedding $(n-1)\times (n-1)$ matrices in the upper left corner of $n\times n$ matrices. Then
$$C_H = \bigcup_{1\leq i\leq n} \big\{ (x_1,\ldots,x_n)\in V^+,\ x_i=0\big\} $$
and the connected components of $V^+\setminus C_H$ are the sets
$$C_i = \big\{ (x_1,\ldots,x_n)\in V^+,\ x_i>0>x_{i+1}\big\} ,$$
where $1\leq i\leq n-1$.
The opposition involution~$\iota : V^+\rightarrow V^+$ is given by
$$\iota(x_1,\ldots,x_n) = (-x_n,\ldots,-x_1)\,;$$
it maps $C_i$ to $C_{n-i}$ for all $1\leq i\leq n-1$.
Here is a reformulation of Theorem~\ref{theoreme principal, general} in the present situation.

\begin{prop}\label{theoreme pour SLn}
Let $\Gamma$ be a discrete subgroup of~$\SLn(\kkk)$ that acts properly discontinuously on $\SLn(\kkk)/\SL_{n-1}(\kkk)$ and that is not a torsion group.
There exists an integer $1\leq i\leq n-1$ such that $\mu(\gamma)\in C_i\cup C_{n-i}$ for almost all~$\gamma\in\Gamma$.
If $\Gamma$ is not virtually cyclic, then $C_i=C_{n-i}$.
\end{prop}

Note that if $n$ is odd, then $C_i\neq C_{n-i}$ for all $1\leq i\leq n-1$, which implies Corollary~\ref{coro SLn}.
Another consequence of Proposition~\ref{theoreme pour SLn} is the following.

\begin{coro}\label{second corollaire SLn}
Assume that $n\geq 4$ is even. Let $\Gamma$ be a discrete subgroup of~$\SLn(\kkk)$ that acts properly discontinuously on $\SLn(\kkk)/\SL_{n-1}(\kkk)$ and that is not virtually cyclic.
Every element $\gamma\in\Gamma$ of infinite order has $n/2$ eigenvalues $t$ with $|t|>1$ and $n/2$ eigenvalues $t$ with $|t|<1$, counting multiplicities.
\end{coro}

The eigenvalues of an element~$g\in\SLn(\kkk)$ belong to some finite extension~$\kkk_g$ of~$\kkk$; in Corollary~\ref{second corollaire SLn} we denote by~$|\cdot|$ the unique absolute value on~$\kkk_g$ extending the absolute value~$|\cdot|$ on~$\kkk$.
As above, replacing~$\kkk$ by~$\kkk_g$, we obtain a Cartan decomposition $\SLn(\kkk_g)=K_gA_g^+K_g$ with $K=K_g\cap\SLn(\kkk)$ and $A^+=A_g^+\cap\nolinebreak\SLn(\kkk)$.
The corresponding Cartan projection $\mu_g : \SLn(\kkk_g)\rightarrow V^+$ extends~$\mu$.

\medskip

\noindent
\textbf{Proof of Corollary~\ref{second corollaire SLn}.}\quad\ 
We may assume that $\Gamma$ is not a torsion group.
Since the only connected component of $V^+\setminus C_H$ that is invariant under~$\iota$ is~$C_{n/2}$, Proposition~\ref{theoreme pour SLn} implies that $\mu(\gamma)\in C_{n/2}$ for almost all $\gamma\in\Gamma$.
Fix an element $\gamma\in\nolinebreak\Gamma$ of infinite order.
Since~$\Gamma$ is discrete and since~$\mu$ is a proper map, $\Vert\mu(\gamma^m)\Vert\rightarrow\nolinebreak +\infty$ as $m\rightarrow +\infty$.
Therefore
$$\frac{1}{m}\,\mu(\gamma^m)\in C_{n/2}$$
for almost all $m\geq 1$.
Let $\lambda:\SLn(\kkk)\rightarrow V^+$ be the \emph{Lyapunov projection} of~$\SLn(\kkk)$, mapping $g\in\SLn(\kkk)$ to~$\mu_g(a_g)$, where $a_g\in\SLn(\kkk_g)$ is any diagonal matrix whose entries are the eigenvalues of~$g$ counted with multiplicities.
By~\cite{be2},~Cor.~2.5,
$$\lambda(\gamma) = \lim_{m\rightarrow +\infty} \frac{1}{m}\,\mu(\gamma^m).$$
Thus $\lambda(\gamma)$ belongs to the closure of~$C_{n/2}$ in~$V^+$.
We claim that $\lambda(\gamma)\notin C_H$.
Indeed, by \cite{be2}, Lem.~4.6, there is a constant $C_{\gamma}>0$ such that for all~$m\geq 1$,
\begin{equation}\label{lambda et mu}
\Vert\lambda(\gamma^m)-\mu(\gamma^m)\Vert \leq C_{\gamma}.
\end{equation}
If $\lambda(\gamma)\in C_H$, then $\lambda(\gamma^m)=m\lambda(\gamma)\in C_H$ for all $m\geq 1$, so that (\ref{lambda et mu}) would contradict the properness criterion (see Section~3.1). This proves the claim.
Therefore, $\lambda(\gamma)\in C_{n/2}$, which means that $\gamma$ has $n/2$ eigenvalues $t$ with $|t|>1$ and $n/2$ eigenvalues $t$ with $|t|<1$, counting multiplicities.
\hfill\qedsymbol

\section{An application to $(G\times G)/\Delta_G$ in the rank-one case}\label{(G*G)/Delta_G}

In this section we prove Theorem~\ref{coro, rang un}, we show that the hypothesis that $\Gamma$ is not a torsion group is necessary in the case of a local field of positive characteristic, and we describe an application to three-dimensional quadrics over a local field.

\subsection{Proof of Theorem~\ref{coro, rang un}}\label{Dem. coro. rang un}
Assume that $\mathrm{rank}_{\kkk}(\mathbf{G})=1$ and let $\mathbf{\Delta_G}$ denote the diagonal of $\mathbf{G}\times\mathbf{G}$.
Fix a Cartan projection~$\mu$ of~$G$ and let $\mu_{\bullet}=\mu\times\mu$ be the corresponding Cartan projection of $G\times G$.
We identify the cone $V^+$ with $\R^+\times\R^+$, and $C_{\Delta_G}$ with the diagonal of $\R^+\times\R^+$.
There are two connected components in $V^+\setminus C_{\Delta_G}$; let $C_+$ (resp.\ $C_-$) denote the one above (resp.\ below) the diagonal.
The opposition involution~$\iota$ is the identity.

We now give a proof of Theorem~\ref{coro, rang un}.
Let $\Gamma$ be a discrete subgroup of~$G\times G$ that acts properly discontinuously on $(G\times G)/\Delta_G$ and that is not a torsion group.
Since $\iota$ is the identity, Theorem~\ref{theoreme principal, general} implies that either $\mu_{\bullet}(\gamma)\in C_+$ for almost all $\gamma\in\Gamma$, or $\mu_{\bullet}(\gamma)\in\nolinebreak C_-$ for almost all $\gamma\in\Gamma$.
Up to switching the factors of~$G\times G$, we may assume that $\mu_{\bullet}(\gamma)\in C_-$ for almost all $\gamma\in\Gamma$.

Let $\mathrm{pr}_1$ (resp.\ $\mathrm{pr_2}$) denote the projection of~$\Gamma$ on the first (resp.\ second) factor of $G\times G$.
The kernel~$F$ of~$\mathrm{pr}_1$ is finite.
If $\Gamma$ is residually finite, then $\Gamma$ contains a normal finite-index subgroup~$\Gamma'$ such that $\Gamma'\cap F$ is trivial.
If~$\Gamma$ is torsion-free, then $F$ is already trivial and we set $\Gamma'=\Gamma$.
In both cases, if we set $\Gamma_0=\mathrm{pr}_1(\Gamma')$, then $\varphi = \mathrm{pr}_2\circ\mathrm{pr}_1^{-1} : \Gamma_0\rightarrow G$ is a group homomorphism and
$$\Gamma' = \{ (g,\varphi(g)),\ g\in\Gamma_0\} .$$
Since $\mu(\varphi(g))<\mu(g)$ for almost all $g\in\Gamma_0$, the group $\Gamma_0$ is discrete in~$G$. Indeed, if it were not, then there would be a sequence $(g_n)_{n\in\N}$ of pairwise distinct points of~$\Gamma_0$ converging to~$1$. Since~$\Gamma$ is discrete in $G\times G$ and since~$\mu$ is a proper map, the sequence $(\mu(\varphi(g_n)))_{n\in\N}$ would tend to infinity. Therefore there would be infinitely many elements $(g,\varphi(g))\in\Gamma$ with $\mu(\varphi(g))\geq\mu(g)$, contradicting the assumption that $\mu_{\bullet}(\gamma)\in C_-$ for almost all $\gamma\in\Gamma$.
This proves that $\Gamma_0$ is discrete in~$G$.
Since $\mu(\varphi(g))<\mu(g)$ for almost all $g\in\Gamma_0$, the properness criterion (see Section~3.1) ensures that for all $R>0$, almost all $g\in\Gamma_0$ satisfy $\mu(\varphi(g))<\mu(g)-R$.

Conversely, if there exist a discrete subgroup $\Gamma_0$ of~$G$ and a group homomorphism $\varphi : \Gamma_0\rightarrow G$ satisfying the conditions of Theorem~\ref{coro, rang un}, then $\Gamma$ acts properly discontinuously on $(G\times G)/\Delta_G$ by the properness criterion.

\subsection{Infinite torsion groups in positive characteristic}\label{Torsion en caracteristique positive}
Take $\mathbf{G}=\mathbf{SL}_2$ over $\kkk=\F_q((t))$, where $\F_q$ is a finite field of characteristic~$p$.
We now give an example of an infinite discrete torsion subgroup of $G\times G$ that acts properly discontinuously on $(G\times G)/\Delta_G$ and nevertheless does not satisfy the conclusions of Theorems~\ref{theoreme principal, general} and~\ref{coro, rang un}.
The Cartan decomposition $G=KA^+K$ holds, where $K=\SL_2(\mathcal{O})=\SL_2(\F_q[[t]])$ and where $A^+$ is the set of diagonal matrices $\mathrm{diag}(a_1,a_2)$ of~$G$ with $|a_1|\geq |a_2|$ (see Section~\ref{Application a SL_n/SL_{n-1}}).
Let $\mu$ be the corresponding Cartan projection.
For every $n\in\N$, set
$$g_n = \begin{pmatrix} 1 & t^{-n} \\ 0 & 1 \end{pmatrix}.$$
Note that for $1\leq r\leq p-1$,
$$\mu(g_n^r)=2n.$$
This can be seen by expanding $g_n^r$ as follows:
$$g_n^r = \begin{pmatrix} 1 & rt^{-n} \\ 0 & 1 \end{pmatrix} = \begin{pmatrix} r & 0 \\ t^n & r^{-1} \end{pmatrix} \begin{pmatrix} t^{-n} & 0 \\ 0 & t^n \end{pmatrix} \begin{pmatrix} r^{-1}t^n & 1 \\ -1 & 0 \end{pmatrix}.$$
Let $\Gamma$ be the subgroup of $G\times G$ generated by the elements $(g_n,g_{2n})$ and the elements $(g_{2n},g_n)$, where $n\in\N$.
It is an infinite residually finite discrete subgroup of~$G$ and each of its nontrivial elements has order~$p$.
The group~$\Gamma$ acts properly discontinuously on $(G\times G)/\Delta_G$ by the properness criterion (see Section~3.1).
It is not virtually cyclic.
However, the two connected components of $V^+\setminus C_{\Delta_G}$ both contain infinitely many points of the form~$\mu(\gamma)$, where~$\gamma\in\Gamma$.

\subsection{An application to three-dimensional quadrics over a local field}
As was pointed out in the introduction, one of the motivations for our investigation of $(G\times G)/\Delta_G$ in the rank-one case is its application to three-dimensional quadrics over a local field~$\kkk$. We now discuss this point in more~detail.

Let $\kkk$ be a local field and $Q$ be a quadratic form on~$\kkk^4$. Consider the unit sphere
$$S(Q) = \{ x\in\kkk^4,\ Q(x)=1\} .$$
By Witt's theorem, it identifies with the homogeneous space $\SO(Q)/H$, where $\mathbf{SO}(Q)$ is the special orthogonal group of~$Q$ and $\mathbf{H}$ is an algebraic $\kkk$-subgroup of~$\mathbf{SO}(Q)$ defined as the stabilizer of some point $x\in S(Q)$.

If $Q$ is $\kkk$-anisotropic, then $\SO(Q)$ is compact (\cite{bot}, \S~4.24); thus every discrete subgroup of~$\SO(Q)$ is finite and acts properly discontinuously on~$S(Q)$.

Assume that $Q$ has Witt index one, \emph{i.e.}, that $\mathrm{rank}_{\kkk}(\mathbf{SO}(Q))=1$.
If $\mathbf{H}$ is $\kkk$-anisotropic, then $H$ is compact, and every discrete subgroup of~$\SO(Q)$ acts properly discontinuously on~$S(Q)$.
On the other hand, if $\mathrm{rank}_{\kkk}(\mathbf{H})=\nolinebreak 1$, then every discrete subgroup of~$\SO(Q)$ acting properly discontinuously on~$S(Q)$ is finite: this is the Calabi-Markus phenomenon (\cite{ko1}, Cor.~4.4).
For instance, if~$\kkk=\R$, then every quadratic form $Q$ on~$\kkk^4$ of Witt index one is equivalent to $x_1^2 - x_2^2 - x_3^2 - x_4^2$ or to $x_1^2 + x_2^2 + x_3^2 - x_4^2$. In the first case, $\SO(Q)$ (resp.~$H$) is isomorphic to~$\SO(1,3)$ (resp.\ to~$\SO(3)$) and every discrete subgroup of~$\SO(Q)$ acts properly discontinuously on~$S(Q)$. In the second case, $\SO(Q)$ (resp.~$H$) is isomorphic to~$\SO(3,1)$ (resp.\ to~$\SO(2,1)$) and every discrete subgroup of~$\SO(Q)$ acting properly discontinuously on~$S(Q)$ is finite.

Now assume that $Q$ has Witt index two, \emph{i.e.}, that $\mathrm{rank}_{\kkk}(\mathbf{SO}(Q))=2$.
For instance, if $\kkk=\R$, then $\SO(Q)$ (resp.\ $H$) is isomorphic to~$\SO(2,2)$ (resp.\ to~$\SO(1,2)$).
We may assume that $Q$ is given by
$$Q(x_1,x_2,x_3,x_4) = x_1 x_4 - x_2 x_3$$
and that~$\mathbf{H}$ is the stabilizer of~$x=(1,0,0,1)\in S(Q)$.
Note that there is a natural transitive action of the group $\SL_2(\kkk)\times\SL_2(\kkk)$ on~$S(Q)$.
Indeed, $\SL_2(\kkk)\times\SL_2(\kkk)$ acts on~$\M_2(\kkk)$ by the formula $(g_1,g_2)\cdot u = g_1 u\, g_2^{-1}$ for all $(g_1,g_2)\in\SL_2(\kkk)\times\SL_2(\kkk)$ and all $u\in\M_2(\kkk)$; identifying $\M_2(\kkk)$ with~$\kkk^4$ gives a linear action of $\SL_2(\kkk)\times\SL_2(\kkk)$ on~$\kkk^4$ that preserves~$Q$ and is transitive on~$S(Q)$.
Since the stabilizer of~$x=(1,0,0,1)$ in $\SL_2(\kkk)\times\nolinebreak\SL_2(\kkk)$ \linebreak is~$\Delta_{\SL_2(\kkk)}$, the quadric~$S(Q)$ identifies with the homogeneous space \linebreak $(\SL_2(\kkk)\times\SL_2(\kkk))/\Delta_{\SL_2(\kkk)}$.
By Theorem~\ref{coro, rang un}, up to switching the factors of $\SL_2(\kkk)\times\SL_2(\kkk)$, the torsion-free discrete subgroups~$\Gamma$ of $\SL_2(\kkk)\times\SL_2(\kkk)$ acting properly discontinuously on $S(Q)$ are exactly the graphs of the form
$$\Gamma = \{ (\gamma,\varphi(\gamma)),\ \gamma\in\Gamma_0\} ,$$
where $\Gamma_0$ is a discrete subgroup of~$\SL_2(\kkk)$ and $\varphi : \Gamma_0\rightarrow\SL_2(\kkk)$ is a group homomorphism such that for all $R>0$, almost all $\gamma\in\Gamma_0$ satisfy $\mu(\varphi(\gamma))<\nolinebreak\mu(\gamma)-\nolinebreak R$.

\vspace{0.5cm}


\begin{thebibliography}{99}

{\normalsize

\bibitem[Alp]{alp}
\textsc{R. C. Alperin}, \textit{An elementary account of Selberg's lemma}, Enseign. Math. (2) 33 (1987), No. 3--4, p. 269--273.

\medskip

\bibitem[Be1]{be1}
\textsc{Y. Benoist}, \textit{Actions propres sur les espaces homog\`enes r\'eductifs}, Ann. Math. 144 (1996), p. 315--347.

\medskip

\bibitem[Be2]{be2}
\textsc{Y. Benoist}, \textit{Propri\'et\'es asymptotiques des groupes lin\'eaires}, Geom. Funct. Anal. 7 (1997), p. 1--47.

\medskip

\bibitem[BoT]{bot}
\textsc{A. Borel, J. Tits}, \textit{Groupes r\'eductifs}, Inst. Hautes \'Etudes Sci. Publ. Math. 27 (1965), p. 55--150.

\medskip

\bibitem[BT1]{bt1}
\textsc{F. Bruhat, J. Tits}, \textit{Groupes r\'eductifs sur un corps local : I. Donn\'ees radicielles valu\'ees}, Inst. Hautes \'Etudes Sci. Publ. Math. 41 (1972), \linebreak p.~5--251.

\medskip

\bibitem[BT2]{bt2}
\textsc{F. Bruhat, J. Tits}, \textit{Groupes r\'eductifs sur un corps local : II. Sch\'emas en groupes. Existence d'une donn\'ee radicielle valu\'ee}, Inst. Hautes \'Etudes Sci. Publ. Math. 60 (1984), p. 5--184.

\medskip

\bibitem[Che]{che}
\textsc{C. Chevalley}, \textit{Th\'eorie des groupes de Lie : II. Groupes alg\'ebriques}, Actualit\'es Sci. Ind. 1152, Hermann \& Cie, Paris, 1951.

\medskip

\bibitem[CR]{cr}
\textsc{C. Curtis, I. Reiner}, \textit{Representation theory of finite groups and associative algebras}, Pure and Applied Mathematics 11, Interscience Publishers, John Wiley \& Sons, New York, London, 1962.

\medskip

\bibitem[Hel]{hel}
\textsc{S. Helgason}, \textit{Differential geometry, Lie groups, and symmetric spaces}, corrected reprint of the 1978 original, Graduate Studies in Mathematics 34, American Mathematical Society, Providence, RI, 2001.

\medskip

\bibitem[Kar]{kar}
\textsc{F. I. Karpelevich}, \textit{Surfaces of transitivity of semisimple groups of motions of a symmetric space}, Soviet Math. Dokl. 93 (1953), p. 401--404.

\medskip

\bibitem[Kli]{kli}
\textsc{B. Klingler}, \textit{Compl\'etude des vari\'et\'es lorentziennes \`a courbure constante}, Math. Ann. 306 (1996), No. 2, p. 353--370.

\medskip

\bibitem[Ko1]{ko1}
\textsc{T. Kobayashi}, \textit{Proper action on a homogeneous space of reductive type}, Math. Ann. 285 (1989), p. 249--263.

\medskip

\bibitem[Ko2]{ko2}
\textsc{T. Kobayashi}, \textit{On discontinuous groups acting on
homogeneous spaces with noncompact isotropy subgroups}, J. Geom. Phys. 12 (1993), No. 2, p. 133--144.

\medskip

\bibitem[Ko3]{ko3}
\textsc{T. Kobayashi}, \textit{Criterion for proper actions on homogeneous spaces of reductive groups}, J. Lie Theory 6 (1996), No. 2, p. 147--163.

\medskip

\bibitem[KR]{kr}
\textsc{R. S. Kulkarni, F. Raymond}, \textit{3-dimensional Lorentz space-forms and Seifert fiber spaces}, J. Differential Geom. 21 (1985), No. 2, \linebreak p.~231--268.

\medskip

\bibitem[Lan]{lan}
\textsc{E. Landvogt}, \textit{Some functorial properties of the Bruhat-Tits building}, J. Reine Angew. Math. 518 (2000), p. 213--241.

\medskip

\bibitem[Mos]{mos}
\textsc{G. D. Mostow}, \textit{Some new decomposition theorems for semi-simple groups}, Mem. Amer. Math. Soc. 14 (1955), p. 31--54.

\medskip

\bibitem[Neu]{neu}
\textsc{J. Neukirch}, \textit{Algebraic number theory}, Grundlehren der Mathematischen Wissenschaften 322, Springer-Verlag, Berlin, 1999.

\medskip

\bibitem[Rou]{rou}
\textsc{G. Rousseau}, \textit{Euclidean buildings}, in \textit{Nonpositive curvature geometry, discrete groups and rigidity}, Proceedings of the 2004 summer school at the Joseph Fourier Institute in Grenoble, S\'eminaires et Congr\`es 18, Soci\'et\'e Math\'ematique de France, Paris, to appear.

\medskip

\bibitem[Sal]{sal}
\textsc{F. Salein}, \textit{Vari\'et\'es anti-de Sitter de dimension 3 exotiques}, Ann. Inst. Fourier 50 (2000), No. 1, p. 257--284.

\medskip

\bibitem[Sel]{sel}
\textsc{A. Selberg}, \textit{On discontinuous groups in higher-dimensional symmetric spaces}, in \textit{Collected papers}, vol.~1, p.~475--492, Springer-Verlag, Berlin,~1989.

\medskip

\bibitem[Tits]{tits}
\textsc{J. Tits}, \textit{Free subgroups in linear groups}, J. Algebra 20 (1972), \linebreak p.~250--270.}

\vspace{0.5cm}

\end{thebibliography}
\end{document}